\documentclass[11pt,twoside]{article}
\usepackage{amssymb}\usepackage{epsfig}

\newcommand{\ignore}[1]{}
\ignore{The text in between these parentheses is ignored by LaTeX}

\def\@begintheorem#1#2{\par\bgroup{\sc #1\ #2. }\it\ignorespaces}
\def\@opargbegintheorem#1#2#3{\par\bgroup{\sc #1\ #2\ (#3). } \it\ignorespaces}
\def\@endtheorem{\egroup}
\newtheorem{theorem}{Theorem}[section]
\newtheorem{corollary}[theorem]{Corollary}
\newtheorem{lemma}[theorem]{Lemma}
\newtheorem{definition}[theorem]{Definition}
\newtheorem{proposition}[theorem]{Proposition}
\newtheorem{example}[theorem]{Example}
\newtheorem{remark}[theorem]{Remark}
\newcommand{\bt}[1]{\begin{theorem}\label{#1}}
\newcommand{\bc}[1]{\begin{corollary}\label{#1}}
\newcommand{\bl}[1]{\begin{lemma}\label{#1}}
\newcommand{\bd}[1]{\begin{definition}\rm\label{#1}}
\newcommand{\bp}[1]{\begin{proposition}\label{#1}}
\newcommand{\be}[1]{\begin{example}\rm\label{#1}}
\newcommand{\br}[1]{\begin{remark}\rm\label{#1}}
\newcommand{\bpr}{\par{\it Proof}. \ignorespaces}
\newcommand{\et}{\end{theorem}}
\newcommand{\ec}{\end{corollary}}
\newcommand{\el}{\end{lemma}}
\newcommand{\ed}{\end{definition}}
\newcommand{\ep}{\end{proposition}}
\newcommand{\ee}{\end{example}}
\newcommand{\er}{\end{remark}}
\newcommand{\epr}{{\ \vbox{\hrule\hbox{%
\vrule height1.3ex\hskip0.8ex\vrule}\hrule
}}\\\par}

\def \Z {\mathbb{Z}}
\def \C {\mathbb{C}}
\def \ideal {\mbox{ideal}}
\def \n {nowhere-zero }

\begin{document}

\title{\bf Nowhere-Zero Flow Polynomials}

\author{Shmuel Onn
\thanks{Supported in part by a grant from ISF - the
Israel Science Foundation, by a VPR grant at the Technion,
by the Fund for the Promotion of Research at the Technion,
and at MSRI by NSF grant DMS-9810361.}}

\date{}
\maketitle

\begin{abstract}
In this article we introduce the {\em flow polynomial} of a digraph and use
it to study \n flows from a commutative algebraic perspective. Using
Hilbert's Nullstellensatz, we establish a relation between \n flows and
dual flows. For planar graphs this gives a relation between \n flows and
flows of their planar duals. It also yields an appealing proof that every
bridgeless triangulated graph has a \n four-flow.
\end{abstract}

\section{Introduction}

The theory of nowhere-zero flows (see \cite{Sey,Wel} for recent surveys)
was introduced by Tutte \cite{Tut} as an extension of
Tait's earlier work \cite{Tai} on the four-color problem for planar graphs.

Let $G=(V,E)$ be a digraph and let $p\geq 2$ be an integer.
A {\em $p$-flow} on $G$ is a mapping $\phi:E\longrightarrow\Z_p$ from arcs
to the additive group $\Z_p=\{0,1,\dots,p-1\}$ of integers modulo $p$
such that preservation holds at each vertex $v$, that is
$\sum\{\phi(e):e\in\delta^-(v)\}-\sum\{\phi(e):e\in\delta^+(v)\}=0$
in $\Z_p$, where $\delta^-(v)\ ,\delta^+(v)$ are the sets of arcs with
head $v$ and tail $v$ respectively. It is a {\em \n $p$-flow} if
$\phi(E)\subseteq \Z_p^*:=\{1,\dots,p-1\}$. If $\phi$ is a flow then
the (signed) sum of arc values on each cocircuit of $G$ is $0$ in $\Z_p$.
With matroid duality in mind, we call $\phi$ a {\em dual $p$-flow}
if the sum of arc values on each circuit of $G$ is zero, and call it
{\em \n dual $p$-flow} if it is \n$\!\!$.

An undirected graph $G$ will be called {\em $p$-flowing} if some orientation
of $G$ admits a \n $p$-flow (and hence so does every orientation - just
flip the sign of $\phi(e)$ whenever $e$ is flipped). Likewise,
$G$ is {\em dually $p$-flowing} if some (and hence every) orientation
of $G$ admits a dual \n $p$-flow. The first fact that motivates flow
theory is the following relation between dual flow and coloring,
implicit in the aforementioned work of Tait \cite{Tai}.
We outline the simple proof.
\bp{Coloring}
A graph is dually $p$-flowing if and only if it is $p$-colorable.
\ep
\bpr
Assume $G$ is connected and oriented with a suitable dual \n $p$-flow $\phi$.
Pick a spanning tree $T$ and a vertex $v$. Set $\omega(v):=0$ and for each
other vertex $u$ set $\omega(u)$ to be the (signed) sum in $\Z_p$ of the values
$\phi(e)$ on arcs on the unique path in $T$ from $v$ to $u$. Since $\phi$ sums
to zero on each circuit, for every arc $e=ab$ we get
$\omega(b)-\omega(a)=\phi(e)$ and since $\phi$ is \n it follows that the
resulting $\omega:V\longrightarrow\Z_p$ is a $p$-coloring.
The converse is likewise easy to see.
\epr
A graph can be flowing only if it has no coloop (also called cut-edge,
isthmus, or bridge); and it can be dually flowing only if it has no loop.
Two gems of flow theory are the following. First, Tutte conjectured
\cite{Tut} that every bridgeless graph is $5$-flowing;
while this is still open, Seymour has shown that every bridgeless graph
is indeed $6$-flowing, see \cite{Sey}. Second, if $G$ is a directed
plane graph then a map $\phi$ is a dual flow precisely when it is a
flow of the plane dual $G^*$; the four-color theorem is thus equivalent
to the statement that every bridgeless planar graph is $4$-flowing.

In this article we introduce the {\em flow polynomial} of a digraph and
use it to study flows from a commutative algebraic perspective.
While the general approach follows the line taken by Lov\'asz in studying
stable sets \cite{Lov} and Alon-Tarsi in studying coloring \cite{AT},
here, inspired by our recent work \cite{BOT}, we take a closer look at
a suitable {\em normal form} of the polynomials that arise.
Using Hilbert's Nullstellensatz, we establish a relation
between \n flows and dual flows. For planar graphs this gives a relation
between \n flows and flows of their planar duals. To state it, we need some
more notation. A map $\phi:E\longrightarrow\Z_p$ is {\em even}
if the number $|\phi^{-1}(p-1)|$ of arcs labelled by the maximal
label $p-1$ is even; otherwise it is {\em odd}. Let
$\psi:E\longrightarrow \Z_p^0:=\{0,\dots,p-2\}$ be a nowhere-(p-1) map.
We say that $\phi:E\longrightarrow\Z_p$ is {\em $\psi$-conformal}
if $\phi(e)\in\{\psi(e),p-1\}$ for every arc $e$.
We establish the following theorem.
\bt{EvenOdd}
A digraph has a \n $p$-flow if and only if it has a nowhere-(p-1) map
$\psi$ such that the number of even $\psi$-conformal dual $p$-flows is
not equal to the number of odd ones.
\et
Since planar duality interchanges circuits and cocircuits,
this gives at once the following corollary.
\bc{Planar}
A plane digraph has a \n $p$-flow if and only if
its plane dual has a nowhere-(p-1) map $\psi$ with number of even
$\psi$-conformal $p$-flows different than that of odd ones.
\ec
Another corollary concerns triangulated graphs: while it can be shown
directly, we find the proof below, which gives a stronger statement
on conformal maps of the zero map, particularly elegant.
\bc{Triangulated}
Any bridgeless triangulated (chordal) graph is $4$-flowing.
\ec
\bpr
We prove by induction on the number of edges the following claim: any
undirected bridgeless triangulated $G=(V,E)$ has an orientation $D$ such
that, for the identically zero map $\psi\equiv0$, the map $\phi=\psi\equiv0$
is the only $\psi$-conformal dual $4$-flow.
If $E$ is empty then the claim is trivially true. Otherwise, pick any
circuit $C\subseteq E$ of size $\leq 3$ in $G$. By induction, the
contraction $G':=G/C$ has an orientation $D'$ satisfying the claim. Extend
$D'$ to an orientation $D$ of $G$ by making $C$ a directed cycle. Consider
any $0$-conformal dual $4$-flow $\phi$ on $D$. Then $\phi(e)\in\{0,3\}$
for all $e$, $\sum_{e\in C}\phi(e)=0$ in $\Z_4$, and $|C|\leq 3$ imply that
$\phi(e)=0$ for all $e\in C$. Now let $\phi'$ be the restriction of $\phi$
to $D'$. Then $\phi'$ is a $0$-conformal dual $4$-flow on $D'$ and hence,
by induction, $\phi'(e)=0$ for all $e\in E\setminus C$. Thus, as claimed,
$\phi\equiv0$, and we are done by Theorem \ref{EvenOdd}.
\epr

\section{The flow polynomial of a digraph}

Fix a digraph $G=(V,E)$ and an integer $p\geq 2$.
Let $x=(x_e:e\in E)$ be a tuple of variables indexed by the arcs of $G$,
and let $\C[x]=\C[x_e:e\in E]$ be the algebra of polynomials with complex
coefficients in these variables. We consider the following polynomial ideal
$$I_E^p\quad:=\quad\ideal\left\{\sum_{i=0}^{p-1}x_e^i\ :\ e\in E\right\}\ ,$$
determined by $p$ and the number of arcs,
and we introduce the following {\em flow polynomial} of $G$,
$$f_G^p\quad:=\quad
\prod_{v\in V}\ \sum_{i=0}^{p-1}\left(\prod_{e\in\delta^-(v)}x_e
\prod_{e\in\delta^+(v)}x_e^{p-1}\right)^i\ .$$
In this section we establish the following statement.
\bt{Membership}
A digraph $G=(V,E)$ has a \n $p$-flow if and only if $f_G^p$ is not in $I_E^p$.
\et
The proof will follow from two properties of the ideal and polynomial which
we establish next. A tuple $a=(a_e:e\in E)$ of complex numbers is a
{\em zero} of $I_E^p$ if $f(a)=0$ for all polynomials $f\in I_E^p$.
Throughout, let $\rho:=\exp({2\pi\sqrt-1\over p})$
denote the primitive $p$-th complex root of unity.
\bp{Zeros}
A tuple $a$ is a zero of $I_E^p$ if and only if
$a_e\in\{\rho^1,\dots,\rho^{p-1}\}$ for all $e\in E$.
Moreover, $I_E^p$ is radical and hence consists precisely
of all polynomials vanishing on its zero set.
\ep
\bpr
The univariate polynomial $f:=\sum_{i=0}^{p-1}z^i$ satisfies
$f\cdot(z-1)=z^p-1=\prod_{i=0}^{p-1}(z-\rho^i)$ and hence its roots are
all $p$-th roots of unity but $\rho^0=1$. Since $I_E^p$ is generated by copies
of $f$, one for each variable $x_e$, the first part of the proposition follows.
Since each such generator has no multiple roots, the ideal is radical.
Therefore, by Hilbert's Nullstellensatz, $I_E^p$ consists precisely of all
polynomials vanishing on its zero set, completing the proof of the proposition.
\epr
The proposition establishes a bijection between \n maps
$\phi:E\longrightarrow\Z_p^*$ and zeros $a=(\rho^{\phi(e)}:e\in E)$ of $I_E^p$.
The \n flows are characterized among such maps $\phi$ by the evaluation
of the flow polynomial on the corresponding zeros $a$, as follows.
\bp{Evaluation}
Consider any map $\phi:E\longrightarrow\Z_p^*$ and let
$a=(\rho^{\phi(e)}:e\in E)$ be the corresponding zero of $I_E^p$.
If $\phi$ is a \n $p$-flow on $G$ then $f_G^p(a)=p^{|V|}$;
otherwise $f_G^p(a)=0$.
\ep
\bpr
Let $s(v):=\sum_{e\in\delta^-(v)}\phi(e)-\sum_{e\in\delta^+(v)}\phi(e)\in\Z_p$
be the flow surplus at vertex $v$. Then
$$f_G^p(a)\ =\
\prod_{v\in V}\sum_{i=0}^{p-1}(\prod_{e\in\delta^-(v)}\rho^{\phi(e)}
\prod_{e\in\delta^+(v)}(\rho^{\phi(e)})^{p-1})^i
\ =\ \prod_{v\in V}\sum_{i=0}^{p-1}(\rho^{s(v)})^i\ .$$
Now, if $s(v)\in\Z_p^*$ then $\sum_{i=0}^{p-1}(\rho^{s(v)})^i=0$
(see proof of Proposition \ref{Zeros}), whereas is $s(v)=0$ then
$\sum_{i=0}^{p-1}(\rho^{s(v)})^i=p$. Since $\phi$ is a flow if and only if
$s(v)=0$ for all $v\in V$, the proof is complete.
\epr

{\em Proof of Theorem \ref{Membership}.}
By Proposition \ref{Zeros}, the polynomial $f_G^p$ is in $I_E^p$
if and only if it vanishes on every zero of $I_E^p$, which, by
Proposition \ref{Evaluation}, holds if and only if $G$ has no \n $p$-flow.
\epr
\br{Two}
Note that the flow polynomial $f_G^p$ has the following very special
appealing property: its evaluations on the zero set of $I_E^p$ assume
{\em only two distinct values}, either $0$ or $p^{|V|}$.
\er
\be{Example}
Let $G=(V,E)$ be a digraph consisting of two vertices $v_1,v_2$ and three
arcs $e_1=e_3=v_1v_2$, $e_2=v_2v_1$, and let $p=3$. The flow polynomial is
{\small
$$f_G^3\!=\! \sum_{i=0}^2(x_1^2x_2x_3^2)^i \cdot \sum_{i=0}^2(x_1x_2^2x_3)^i
\!=\!x_1^6x_2^6x_3^6+x_1^5x_2^4x_3^5+x_1^4x_2^5x_3^4+x_1^4x_2^2x_3^4+
x_1^3x_2^3x_3^3+x_1^2x_2^4x_3^2+x_1^2x_2x_3^2+x_1x_2^2x_3+1;$$}$\!$
the zeros of the ideal $I_E^3=\ideal\{x_1^2+x_1+1,x_2^2+x_2+1,x_3^2+x_3+1\}$
are all $8$ tuples $a=(a_1,a_2,a_3)$ with each $a_i\in\{\rho,\rho^2\}$ where
$\rho=\exp({2\pi\sqrt-1\over 3})$; the evaluation of $f_G^3$ on a zero of
$I_E^3$ corresponding to a \n map $\phi=(\phi_1,\phi_2,\phi_3)$ is $3^2=9$
if $\phi$ is either $(1,2,1)$ or $(2,1,2)$, and is $0$ otherwise,
distinguishing $(1,2,1)$ and $(2,1,2)$ as the only two \n $3$-flows of $G$.
\ee

\section{The normal form of the flow polynomial}

We proceed to take a close look at the {\em normal form} of the flow
polynomial with respect to a natural monomial basis of the quotient
$\C[x]/I_E^p$. Consider the following set of {\em basic monomials},
$$B\quad:=\quad \left\{\prod_{e\in E} x_e^{\psi(e)}\quad :\quad
\,\psi:E\longrightarrow\Z_p^0=\{0,\dots,p-2\}\right\}\quad .$$
\bp{Basis}
The (congruence classes of) basic monomials
form a $\C$-basis of $\C[x]/I_E^p$.
\ep
\bpr
First, it is clear that $I_E^p$ contains no nonzero polynomial which is a linear combination
of monomials in $B$, so $B$ is linearly independent modulo $I_E^p$. Second,
$I_E^p$ is radical by Proposition \ref{Zeros}, and hence, by Hilbert's
Nullstellensatz, the vector space dimension of $\C[x]/I_E^p$ equals the
number of zeros of $I_E^p$. By Proposition \ref{Zeros} this number is
$(p-1)^{|E|}$, which is precisely the number of basic monomials, and so
it follows that $B$ spans the quotient space and hence provides a basis.
\epr

It follows that for every polynomial $f\in\C[x]$ there is a unique
polynomial $[f]$, called the {\em normal form} of $f$, which satisfies
$f-[f]\in I_E^p$ and is a $\C$-linear combination of basic monomials
$$[f]\ \ =\ \sum_{\psi:E\rightarrow\Z_p^0}
c_{\psi}\prod_{e\in E} x_e^{\psi(e)}\ \ .$$
In particular, $f\in I_E^p$ if and only if $[f]=0$.
By characterizing the normal form of the flow polynomial we will be able,
via Theorem \ref{Membership}, to establish the promised criterion of
Theorem \ref{EvenOdd} for a graph to be flowing. We proceed to
study normal forms, starting with powers of variables.
\bp{NormalForm}
The normal form of $x_e^{q\cdot p+r}$
with $q$ any nonnegative integer and $r\in\Z_p$ is
$$[x_e^{q\cdot p+r}]=\left\{
\begin{array}{cll}
  x_e^r & \mbox{if} & r\in\Z_p^0\\
  -\sum_{i=0}^{p-2}x_e^i & \mbox{if} & r=p-1
\end{array}\right.$$
\ep
\bpr
First, we have $x_e^p-1=(x_e-1)\cdot\sum_{i=0}^{p-1}x_e^i\in I_E^p$
and $1$ is a basic monomial, which shows that $[x_e^p]=1$.
Thus, the normal form
of an arbitrary power of $x_e$ is determined by the normal form of powers
$x_e^r$ with $r\in\Z_p$. If $r\in\Z_p^0$ then the monomial $x_e^r$ is basic
and hence satisfies $[x_e^r]=x_e^r$. If $r=p-1$ then
$x_e^{p-1}-(-\sum_{i=0}^{p-2}x_e^i)=\sum_{i=0}^{p-1}x_e^i \in I_E^p$
and $-\sum_{i=0}^{p-2}x_e^i$ is a linear combination of basic monomials,
so $[x_e^{p-1}]=-\sum_{i=0}^{p-2}x_e^i$. This completes the proof.
\epr
Now, for any two polynomials $f,g$ and scalars $s,t\in\C$ we have
$[sf+tg]=s[f]+t[g]$ and $[fg]=[[f][g]]$. The first identity implies that the
normal form of any polynomial is determined by the normal forms of its
monomials. The second identity implies that for any monomial
$\prod_{e\in E} x^{m_e}$ we have
$[\prod_{e\in E} x^{m_e}]=[\prod_{e\in E} [x^{m_e}]]$; but
Proposition \ref{NormalForm} implies that the polynomial
$\prod_{e\in E} [x^{m_e}]$ is in the $\C$-linear span of basic monomials and
hence $[\prod_{e\in E} x^{m_e}]=\prod_{e\in E} [x^{m_e}]$. This completely
determines the normal form of any polynomial.

\vspace*{.5cm}\noindent{\bf Example \ref{Example} continued.}
Consider again the digraph $G$ with three arcs, and let again $p=3$.
Using Proposition \ref{NormalForm} we find that the
normal form of the flow polynomial is
{\small\begin{eqnarray*}
[f_G^3]&=& 1+(-x_1-1)x_2(-x_3-1)+x_1(-x_2-1)x_3+x_1(-x_2-1)x_3+1\\
&+&(-x_1-1)x_2(-x_3-1)+(-x_1-1)x_2(-x_3-1)+x_1(-x_2-1)x_3+1\\
&=&3(x_1x_2-x_1x_3+x_2x_3+x_2+1)\ .
\end{eqnarray*} }
$\!\!$Since $[f_G^3]\neq 0$, we find that $f_G^3\notin I_E^3$ and hence,
by Theorem \ref{Membership}, $G$ admits a \n $3$-flow.

\vspace*{.6cm}
We next show that the coefficients of the monomials in the normal form of the
flow polynomial can be nicely interpreted in terms of certain dual flows.
Recall that a map $\phi:E\longrightarrow\Z_p$ is {\em even} if the number
$|\phi^{-1}(p-1)|$ of arcs labelled by the maximal label $p-1$ is even;
otherwise it is {\em odd}. Recall also that $\phi$ is {\em $\psi$-conformal}
for a nowhere-(p-1) map $\psi:E\longrightarrow \Z_p^0=\{0,\dots,p-2\}$ if
$\phi(e)\in\{\psi(e),p-1\}$ for every arc $e$. We have the following theorem.
\bt{Main}
Let $G=(V,E)$ be an orientation of a connected graph,
and let $p\geq 2$ be an integer.
Then the normal form of the flow polynomial of $G$ is given by
$$[f_G^p]\ \ =\ p\cdot\sum_{\psi:E\longrightarrow\Z_p^0}
c(\psi)\prod_{e\in E} x_e^{\psi(e)}\ \ ,$$
where $c(\psi)$ denotes the number of even $\psi$-conformal
dual $p$-flows minus the number of odd ones.
\et
\bpr
The flow polynomial can be expanded as
$$f_G^p\quad:=\quad \sum_{\omega:V\longrightarrow\Z_p}
\ \prod_{v\in V}\ \left(\prod_{e\in\delta^-(v)}x_e
\prod_{e\in\delta^+(v)}x_e^{p-1}\right)^{\omega(v)}\ ,$$
the sum extending over all labellings $\omega$ of vertices by $\{0,\dots,p-1\}$.
Since each arc $e=uv$ satisfies $e\in\delta^+(u)$ and $e\in\delta^-(v)$
we can rewrite this as
$$f_G^p\quad:=\quad \sum_{\omega:V\longrightarrow\Z_p}
\ \prod_{e=uv\in E}\ x_e^{\omega(v)}x_e^{(p-1)\omega(u)}\ .$$
Consider any $\omega:V\longrightarrow\Z_p$ and let
$\phi:E\longrightarrow\Z_p$ be the map that labels each arc
$e=uv$ by $\phi(e):=\omega(v)-\omega(u)$ in $\Z_p$.
By Proposition \ref{NormalForm}, we then have
$[x_e^{\omega(v)}x_e^{(p-1)\omega(u)}]=[x_e^{\phi(e)}]$
and hence the normal form of the summand in the above expression
of $f_G^p$ corresponding to $\omega$ satisfies
$$\left[\prod_{e=uv\in E}\ x_e^{\omega(v)}x_e^{(p-1)\omega(u)}\right]
\quad=\quad \prod_{e\in E}\left[x_e^{\phi(e)}\right]\ .$$
Now, since the arc labelling $\phi$ is induced from a vertex
labelling $\omega$, the (signed) sum of the $\phi$ values of arcs on each
circuit of $G$ is $0$ in $\Z_p$ and hence $\phi$ is a dual $p$-flow.
Since the undirected graph underlying $G$ is connected, $\omega$ is
uniquely determined by $\phi$ and the value $\omega(v)$ on an arbitrary
vertex $v$ (see proof of Proposition \ref{Coloring}), so $\phi$
arises from precisely $p$ distinct maps $\omega$, and we get
\begin{eqnarray*}
[f_G^p]& = & p\cdot \sum\left\{\prod_{e\in E}\left[x_e^{\phi(e)}\right]
\quad :\quad \mbox{$\phi$ dual $p$-flow}\right\}\\
& =& p\cdot \sum\left\{\prod_{e:\phi(e)\in\Z_p^0} x_e^{\phi(e)}
\prod_{e:\phi(e)=p-1}(-\sum_{i\in\Z_p^0} x_e^i)
\quad :\quad \mbox{$\phi$ dual $p$-flow}\right\}\ \ .
\end{eqnarray*}
Now consider the basic monomial $\prod_{e\in E} x_e^{\psi(e)}$
corresponding to $\psi:E\longrightarrow \Z_p^0$. Then, in the right
hand side sum in the above expression of $[f_G^p]$,
every even $\psi$-conformal dual $p$-flow map $\phi$ contributes a term
$\prod_{e\in E} x_e^{\psi(e)}$, whereas every odd one contributes a term
$-\prod_{e\in E} x_e^{\psi(e)}$. This shows that, as claimed, the
coefficient $c(\psi)$ of $\prod_{e\in E} x_e^{\psi(e)}$ in $[f_G^p]$
is equal to the number of even $\psi$-conformal dual $p$-flows minus
the number of odd ones, completing the proof of the theorem.
\epr
\br{Components}
More generally, if $G$ is an orientation of a graph with $\kappa$ connected
components then a suitable adjustment of the analysis above shows that the
normal form of the flow polynomial is
$$[f_G^p]\ \ =\ p^\kappa \cdot\sum_{\psi:E\longrightarrow\Z_p^0}
c(\psi)\prod_{e\in E} x_e^{\psi(e)}\ \ .$$
\er

\vspace*{.2cm}\noindent{\bf Example \ref{Example} continued.}
Consider once again the digraph with three arcs and let $p=3$.
Let us examine some of the monomials
$x_1^{\psi_1}x_2^{\psi_2}x_3^{\psi_3}$ of the normal form of the
flow polynomial described explicitly before and see how
they can be computed using Theorem \ref{Main}.
For instance, for $\psi=(1,1,0)$, the only $\psi$-conformal dual $3$-flow
is the even $\phi=(2,1,2)$, so the coefficient of $x_1x_2$ in $[f_G^3]$
is $3\cdot 1=3$; for $\psi=(1,0,1)$, the only conformal dual flow is
the odd $(2,1,2)$, so the coefficient of $x_1x_3$ in $[f_G^3]$ is $3\cdot(-1)=-3$;
for $\psi=(1,0,0)$ there is no conformal dual flow so the coefficient of
$x_1$ in $[f_G^3]$ is $0$; for $\psi=(0,0,0)$, the only conformal dual flow is the
even $(0,0,0)$ so the coefficient of $1$ in $[f_G^3]$ is $3$;
and similarly for the four other remaining monomials.

Now consider the plane dual $G^*=(U,E^*)$ of $G$ under some plane embedding
of $G$; it has three vertices $u_1,u_2,u_3$ and three dual arcs
$e_1^*=u_1u_2, e_2^*=u_1u_3, e_3^*=u_2u_3$. The normal form $[f_{G^*}^3]$
of the flow polynomial of $G^*$ can be computed via Theorem \ref{Main}
by considering dual flows of $G^*$, namely flows of $G$.
For instance, for $\psi=(1,1,0)$, the $\psi$-conformal $3$-flows
$\phi$ of $G$ are $(1,1,0),(2,2,0),(2,1,2)$, with all three even, so the
coefficient of $x_1x_2$ in $[f_{G^*}^3]$ is $3\cdot 3=9$; for $\psi=(1,0,1)$,
the conformal flows are $(2,0,1),(1,2,1),(1,0,2)$ with all three odd,
so the coefficient of $x_1x_3$ in $[f_{G^*}^3]$ is $3\cdot(-3)=-9$;
for $\psi=(1,0,0)$, the conformal flows are $(1,0,2),(2,2,0)$ with one odd
and one even, so the coefficient of $x_1$ in $[f_{G^*}^3]$ is $0$;
for $\psi=(0,0,0)$, the conformal flows are $(0,0,0),(2,2,0),(0,2,2)$
with all three even, so the coefficient of $1$ in $[f_{G^*}^3]$ is $9$;
and similarly for the four other remaining monomials, giving
{\small$$
[f_{G^*}^3]\quad=\quad
3(3x_1x_2-3x_1x_3+3x_2x_3+3x_2+3)\quad=\quad3[f_G^3]\ .
$$}
$\!$
\vspace*{.2cm}

We are finally in position to prove the following
theorem stated in the introduction.

\vspace*{.2cm}
\noindent
{\bf Theorem \ref{EvenOdd}}
{\em A digraph has a \n $p$-flow if and only if it has a nowhere-(p-1) map
$\psi$ such that the number of even $\psi$-conformal dual $p$-flows is
not equal to the number of odd ones.}

\vspace*{.2cm}
\bpr
Let $f_G^p$ be the flow polynomial of a digraph $G$.
By Theorem \ref{Membership}, $G$ has a \n $p$-flow if and only if
$f_G^p\notin I_E^p$, which holds if and only if the normal form
$[f_G^p]$ is nonzero. By Theorem \ref{Main}, $[f_G^p]\neq 0$
if and only if $G$ admits a nowhere-(p-1) map $\psi$ such that $c(\psi)\neq 0$.
Since $c(\psi)$ is the number of even $\psi$-conformal dual $p$-flows minus
the number of odd ones, we are done.
\epr

\section{The four-flow polynomial of an undirected graph}

In this section we work out a variant of the flow polynomial for four-flows
for an undirected graph $G=(V,E)$. It is simpler and perhaps better suited
for the study of four-flows of planar graphs and the four-color theorem.
The outline is similar to that of the previous two sections; we therefore
do not go through the proofs which are analogous to those provided before.

Let $G=(V,E)$ be a graph. A {\em four-flow} on $G$ is a mapping
$$\phi=(\phi_1,\phi_2):
E\longrightarrow \Z_2\times\Z_2=\{(0,0),(0,1),(1,0),(1,1)\}$$
such that $\sum\{\phi(e):e\in\delta(v)\}=(0,0)$ for each vertex $v$,
where $\delta(v)$ is the set of edges incident on $v$. As before,
if $\phi$ is a flow then the sum of edge values on each cocircuit of $G$
is $(0,0)$. Again, we call $\phi$ a {\em dual four-flow} if the sum of edge
values on each circuit of $G$ is zero in $\Z_2\times\Z_2$.

Let $x=(x_e:e\in E)$, $y=(y_e:e\in E)$ be two tuples of variables
indexed by edges and let $\C[x,y]$ be corresponding polynomial algebra.
Consider the following ideal and polynomial,
$$I_E:=\ideal\{x_e^2-1,\, y_e^2-1,\, (x_e+1)(y_e+1)\,:\,e\in E\}\ ,$$
$$f_G\quad:=\quad
\prod_{v\in V}\ \left(\prod_{e\in\delta(v)}x_e+1\right)
\left(\prod_{e\in\delta(v)}y_e+1\right)\ .$$
We have the following analog of Theorem \ref{Membership}.
\bt{Membership1}
A graph $G=(V,E)$ has a \n four-flow if and only if $f_G$ is not in $I_E$.
\et
As before , this is a consequence of the following two properties of
$I_E$ and $f_G$. A pair of tuples $a=(a_e:e\in E)$, $b=(b_e:e\in E)$ of
complex numbers is a {\em zero} of $I_E$ if $f(a,b)=0$ for all $f\in I_E$.
\bp{Zeros1}
The pair $(a,b)$ is a zero of $I_E$ if and only if
$(a_e,b_e)\in\{(1,-1),(-1,1),(-1,-1)\}$ for all $e$.
Moreover, $I_E$ is radical and hence consists of all
polynomials vanishing on its zero set.
\ep
The proposition establishes a bijection between \n maps
$$\phi=(\phi_1,\phi_2):
E\longrightarrow(\Z_2\times\Z_2)^*=\{(0,1),(1,0),(1,1)\}$$
and zeros $(a,b)$ of $I_E$ given by
$a_e=(-1)^{\phi_1(e)}$, $b_e=(-1)^{\phi_2(e)}$ for all $e\in E$.
The \n flows are characterized among such maps $\phi$ by the evaluation
of the flow polynomial on the corresponding zeros $(a,b)$, as follows.
\bp{Evaluation1}
Consider any \n map $\phi=(\phi_1,\phi_2)$ and let $(a,b)$ be
the corresponding zero of $I_E$. If $\phi$ is a \n four-flow
on $G$ then $f_G(a,b)=4^{|V|}$; otherwise $f_G(a,b)=0$.
\ep

\vskip.5cm
{\em Proof of Theorem \ref{Membership1}.}
By Proposition \ref{Zeros1}, $f_G$ lies in $I_E$ if and only
if it vanishes on its set of zeros, which, by Proposition \ref{Evaluation1},
holds if and only if $G$ has no \n four-flow.
\epr

As before, we next consider the normal form of the flow polynomial
with respect to a natural monomial basis of the quotient
$\C[x,y]/I_E$. Consider the following set of basic monomials,
$$B\quad:=\quad \left\{\prod_{e\in E} x_e^{\psi_1(e)}y_e^{\psi_2(e)}\quad
:\quad \,
\psi=(\psi_1,\psi_2):E\longrightarrow(\Z_2\times\Z_2)^0:=\{(0,0),(0,1),(1,0)\}
\right\}\quad .$$
\bp{Basis1}
The (congruence classes of) basic monomials
form a $\C$-basis of $\C[x,y]/I_E$.
\ep

For every polynomial $f\in\C[x,y]$ let again $[f]$ denote its normal form
which is the unique $\C$-linear combination of basic monomials satisfying
$f-[f]\in I_E$. The normal form of powers of pairs of variables $x_e,y_e$
is determined by the following analog of Proposition \ref{NormalForm}.
\bp{NormalForm1}
For any two nonnegative integers $q_1,q_2$ and any $r_1,r_2\in\Z_2$ we have
$$[x_e^{2q_1+r_1}y_e^{2q_2+r_2}]=\left\{
\begin{array}{cll}
x_e^{r_1}y_e^{r_2} & \mbox{if} & (r_1,r_2)\in(\Z_2\times\Z_2)^0\\
-x_e-y_e-1 & \mbox{if} & (r_1,r_2)=(1,1)
\end{array}\right.$$
\ep

Now, for any monomial
$\prod_{e\in E} x^{m_e}y^{n_e}$ we have
$[\prod_{e\in E} x^{m_e}y^{n_e}]=[\prod_{e\in E}[x^{m_e}y^{n_e}]]$;
but Proposition \ref{NormalForm1} implies that the polynomial
$\prod_{e\in E} [x^{m_e}y^{n_e}]$ is in the $\C$-linear span of basic monomials
and hence $[\prod_{e\in E} x^{m_e}y^{n_e}]=\prod_{e\in E} [x^{m_e}y^{n_e}]$.
This completely determines the normal form of any monomial and hence,
as explained before, of every polynomial.

\vspace*{.5cm}
We next show, in analogy with Theorem \ref{Main}, an interpretation of the
coefficients of the monomials in the normal form of the flow polynomial
in terms of suitable conformal dual flows. A map $\phi:E\longrightarrow\Z_p$
is even if the number $|\phi^{-1}(1,1)|$ of edges labelled by $(1,1)$
is even; otherwise it is odd. The map $\phi$ is $\psi$-conformal
for a nowhere-(1,1) map
$\psi=(\psi_1,\psi_2):E\longrightarrow(\Z_2\times\Z_2)^0$
if $\phi(e)\in\{\psi(e),(1,1)\}$ for every edge $e$.
We have the following analog of Theorem \ref{Main}.
\bt{Main1}
Let $G=(V,E)$ be a graph with $\kappa$ connected components.
Then the normal form of the four-flow polynomial of $G$ is given by
$$[f_G]\ \ =\ 4^\kappa \cdot
\sum_{\psi=(\psi_1,\psi_2):E\longrightarrow(\Z_2\times\Z_2)^0}
c(\psi)\prod_{e\in E} x_e^{\psi_1(e)}y_e^{\psi_2(e)}\ \ ,$$
where $c(\psi)$ is the number of even $\psi$-conformal
dual four-flows minus the number of odd ones.
\et

\vspace*{.5cm}
We also conclude the following analog of Theorem \ref{EvenOdd}.
\bt{EvenOdd1}
A graph has a \n four-flow if and only if it has a nowhere-(1,1) map
$\psi$ such that the number of even $\psi$-conformal dual four-flows
is not equal to the number of odd ones.
\et
\bpr
Let $f_G$ be the flow polynomial of a graph $G$.
By Theorem \ref{Membership1}, $G$ has a \n four-flow if and only if
$f_G\notin I_E$, which holds if and only if $[f_G]$ is nonzero.
By Theorem \ref{Main1}, $[f_G]\neq 0$
if and only if $G$ admits a nowhere-(1,1) map $\psi$ such that $c(\psi)\neq 0$.
Since $c(\psi)$ is the number of even $\psi$-conformal dual four-flows minus
the number of odd ones, we are done.
\epr

\vskip.2cm\noindent {\small Shmuel Onn}\newline
\emph{Technion - Israel Institute of Technology, 32000 Haifa, Israel.}

\emph{email: onn{\small @}ie.technion.ac.il}

\emph{http://ie.technion.ac.il/{\small $\sim$onn}}

\end{document}